# On quasicomplete $k$-surfaces in $3$-dimensional space-forms

27th November 2022

Graham Smith

**Abstract:** In the study of immersed surfaces of constant positive extrinsic curvature in space-forms, it is natural to substitute completeness for a weaker property, which we here call quasicompleteness. We determine the global geometry of such surfaces under the hypotheses of quasicompleteness. In particular, we show that, for $k > \text{Max}(0, -c)$, the only *quasicomplete* immersed surfaces of constant extrinsic curvature equal to $k$ in the 3-dimensional space-form of constant sectional curvature equal to $c$ are the geodesic spheres. Together with earlier work of the author, this completes the classification of quasicomplete immersed surfaces of constant positive extrinsic curvature in 3-dimensional space-forms.

**1 - Introduction.** When working with immersed surfaces in space-forms, it is natural to impose the condition of completeness. However, when studying surfaces of constant positive extrinsic curvature, various phenomena (see, for example, [8]) indicate that we should instead use the weaker condition of quasicompleteness, which we define as follows. Let $(S, e)$ be an immersed surface in some 3-dimensional space-form $X$, and let $\text{I}_e$, $\text{II}_e$ and $\text{III}_e$ denote its three fundamental forms. We say that $(S, e)$ is *quasicomplete* whenever the riemannian metric $\text{I}_e + \text{III}_e$ is complete. Note that an analogous concept is used in [3], [4], [5] and [7] to study what the authors there call flat fronts.

When studying quasicomplete immersed surfaces of constant extrinsic curvature, it is useful to know what can still be said about the global geometry of their first fundamental forms. It is the purpose of this paper to address this problem. Given $c \in \mathbb{R}$ and $m \in \mathbb{N}$, let $X_c^m$ denote the $m$-dimensional space-form of constant sectional curvature equal to $c$. Choose $k > 0$, and let $(S, e)$ be a quasicomplete immersed surface in $X_c^3$ of constant extrinsic curvature equal to $k$. By Gauss' Theorem, $(S, \text{I}_e)$ has constant curvature equal to $(c + k)$, and thus carries an $X_{c+k}^2$-structure in the sense of Weyl.

We construct the completion of $S$ as follows. We define the pseudodistance $d$ over the set of Cauchy sequences in $S$ by

$$d((x_m)_{m \in \mathbb{N}}, (y_m)_{m \in \mathbb{N}}) := \underset{m \to \infty}{\text{Lim}}\, d(x_m, y_m). \tag{0.1}$$

We identify two Cauchy sequences whenever the pseudodistance between them vanishes. We then define the *completion* $\overline{S}$ to be the space of Cauchy sequences furnished with the metric that $d$ defines. Since every point of $S$ trivially identifies with the constant sequence at that point, $S$ naturally embeds as an open subset of $\overline{S}$, and we denote $\partial S := \overline{S} \setminus S$.

**Theorem 1.1**

*The $X_{c+k}^2$-structure of $S$ extends uniquely to an $X_{c+k}^2$-structure of $\overline{S}$ with geodesic boundary. Furthermore, $e$ extends to a continuous function $\overline{e} : \overline{S} \to X_c^3$ which sends components of $\partial S$ locally isometrically into geodesics in $X_c^3$.*

**Remark 1.1.** We draw the reader's attention to the intriguing analogy between Theorem 1.1 and Corollary E of [1].

In terms of the global geometry of $(S, \text{I}_e)$, this has the following useful consequence.

**Corollary 1.2**

*Any two points of $(S, \text{I}_e)$ are joined by a length-minimising geodesic. Furthermore, when $S$ is simply-connected and $c + k < 0$, this geodesic is unique.*

A nice application of Theorem 1.1 is the following extension of a classical result of Liebmann.

**Theorem 1.3**

When $k > \text{Max}(0, -c)$, the only quasicomplete immersed surfaces in $X_c^3$ of constant extrinsic curvature equal to $k$ are the geodesic spheres.

**Remark 1.2.** The classical version of this result, which holds for complete immersed surfaces, is proven in Theorem III.5.2 of [10] and Theorem V.1.5 of [2].

In Theorem 2.8.3 of [9] we extended the theorem of Volkov–Vladimirova and Sasaki to the quasicomplete case, proving that the only quasicomplete immersed surfaces in $\mathbb{H}^3$ of constant extrinsic curvature equal to $1$ are the horospheres and the level sets of distance functions to complete geodesics. Likewise, in Theorems 1.3.1 and 1.3.2 of [8], we classified quasicomplete immersed surfaces in $\mathbb{H}^3$ of constant extrinsic curvature equal to $0 < k < 1$. Together with Theorem 1.3, this yields a complete classification of quasicomplete immersed surfaces of constant positive extrinsic curvature in 3-dimensional space forms. Significantly, the corresponding classification of *complete* surfaces is unfinished, precisely because of the difficulties arising from the case of surfaces of constant extrinsic curvature equal to $0 < k < 1$ in $\mathbb{H}^3$. We view this as further evidence of the naturality of the quasicompleteness condition.

**2 - Acknowledgements.** This work was carried out during a research visit to the Institut des Hautes Études Scientifiques. The author is grateful for the excellent conditions provided during that stay. The author is grateful to François Labourie for his interest in this work and to Andrea Seppi for helpful comments made to an earlier draft of this paper.

**3 - Labourie space.** For $k > 0$, we define a *$k$-surface* in $X_c^3$ to be a quasicomplete, immersed surface $(S, e)$ of constant extrinsic curvature equal to $k$. Let $S_1 X_c^3$ denote the unit sphere bundle over $X_c^3$. Given a $k$-surface $(S, e)$, we denote by $\hat{e} : S \to S_1 X_c^3$ its outward-pointing, unit, normal vector field and we call $(S, \hat{e})$ its *Gauss lift*. We define a *Gauss lift* in $S_1 X_c^3$ to be any immersed surface obtained from some $k$-surface in this manner. Note that, by definition, every Gauss lift is complete.

Theorem 1.1 will follow in a straightforward manner from properties of the space of Gauss lifts which we now describe. We first recall the concept of Cheeger-Gromov convergence, which serves to define its topology. Let $S_1 X_c^3$ denote the unit sphere bundle over $X_c^3$. We define a *marked surface* in $S_1 X_c^3$ to be a triple $(S, f, p)$, where $(S, f)$ is a complete immersed surface, and $p \in S$. We identify marked surfaces which are reparametrisations of one another. Let $(S_m, f_m, p_m)$ be a sequence of marked surfaces in $S_1 X_c^3$. We say that this sequence *converges* towards the marked surface $(S_\infty, f_\infty, e_\infty)$ in the Cheeger–Gromov sense whenever there exists a sequence $(\Phi_m)$ of functions such that

(1) for all $m$, $\Phi_m : S_\infty \to S_m$ and $\Phi_m(p_\infty) = p_m$; and

for every relatively compact open subset $\Omega$ of $S_\infty$, there exists $M$ such that

(2) for all $m \geqslant M$, $\Phi_m$ restricts to a diffeomorphism of $\Omega$ onto its image; and

(3) $(f_m \circ \Phi_m)_{m \geqslant M}$ converges to $f_\infty$ in the $C^\infty_{\text{loc}}$ sense over $\Omega$.

We call $(\Phi_m)$ a sequence of *convergence maps* of $(S_m, f_m, p_m)$ with respect to $(S_\infty, f_\infty, p_\infty)$.

We define a *marked Gauss lift* in $S_1 X_c^3$ to be a triple $(S, \hat{e}, p)$ where $(S, \hat{e})$ is a Gauss lift and $p \in S$. We define a *marked tube* in $S^1 X_c^3$ to be a triple $(S, \hat{e}, p)$ where $\hat{e} : S \to S_1 X_c^3$ is a covering map of the unit normal bundle over some complete geodesic in $X_c^3$. Let $\text{MGL}_{k,c}$ and $\text{MT}_c$ denote respectively the sets of reparametrisation equivalence classes of marked Gauss lifts and marked tubes in $S_1 X_c^3$. We denote

$$\overline{\text{MGL}}_{k,c} := \text{MGL}_{k,c} \cup \text{MT}_c, \tag{0.2}$$

and we furnish this set with the Cheeger–Gromov topology. $\text{Isom}(X_c^3)$ acts continuously on this space by post-composition, and we define *Labourie space* by

$$\mathcal{L}_{k,c} := \overline{\text{MGL}}_{k,c} / \text{Isom}(X_c^3). \tag{0.3}$$

In [6], Labourie proves the following result.

**Theorem 3.1, Labourie (1997)**

$\mathcal{L}_{k,c}$ is compact.

**Remark 3.1.** The reader may consult our recent review [9] for an alternative proof of Theorem 3.1.

Define $d : \mathrm{MGL}_{k,c} \to ]0, \infty]$ by
$$d(S, \hat{e}, p) := d(p, \partial S), \tag{0.4}$$

where distance is taken with respect to the metric induced by $(\pi \circ \hat{e})$. Note that this function is infinite if and only if $(S, \pi \circ e)$ is complete. We extend $d$ to a function defined over $\overline{\mathrm{MGL}}_{k,c}$ by setting it equal to zero over $\mathrm{MT}_c$. This function is trivially invariant under the action of $\mathrm{Isom}(X_c^3)$, and thus descends to a function over $\mathcal{L}_{k,c}$.

**Lemma 3.2**

$d$ is lower-semicontinuous, that is, for all $x \in \mathcal{L}_{k,c}$,
$$\mathop{\mathrm{LimInf}}_{y \to x} d(y) \geqslant d(x). \tag{0.5}$$

**Remark 3.2.** Note that equidistant surfaces to some fixed geodesic $\Gamma$ in $\mathbb{H}^3$ have constant sectional curvature equal to 1. For any marked Gauss lift $(S, \hat{e}, p)$ of any such surface, $d(S, \hat{e}, p) = \infty$. However, upon letting the radius tend to zero, we obtain sequences converging towards marked tubes. This shows that $d$ is not continuous over $\mathcal{L}_{-1,1}$. Using rotationally symmetric surfaces about geodesics in $\mathbb{H}^3$ (see, for example, Chapter 7.F of [11]), we see that $d$ is likewise not continuous over $\mathcal{L}_{-1,k}$ for $0 < k < 1$.

**Proof:** Let $(S_m, \hat{e}_m, p_m)_{m \in \mathbb{N}}$ be a sequence in $\overline{\mathrm{MGL}}_{k,c}$ converging towards the limit $(S_\infty, \hat{e}_\infty, p_\infty)$ and let $(\Phi_m)$ be a corresponding sequence of convergence maps. We may suppose that $(S_\infty, \hat{e}_\infty, p_\infty)$ is not a tube, for otherwise the result holds trivially. In particular, $r_\infty := d(S_\infty, \hat{e}_\infty, p_\infty) > 0$. For all $r > 0$, let $B_r$ denote the open ball of radius $r$ about $p_\infty$ in $S_\infty$ with respect to the metric induced by $(\pi \circ \hat{e}_\infty)$. Choose $\epsilon > 0$ and note that $B := \overline{B}_{r_\infty - \epsilon}$ is compact. Let $M$ be such that, for all $m \geqslant M$, the restriction of $\Phi_m$ to $B$ is a diffeomorphism onto its image. For all $m \geqslant M$, and for all $r$, let $B_{r,m}$ denote the open ball of radius $r$ about $p_\infty$ with respect to the metric induced by $(\pi \circ e_m \circ \Phi_m)$. For sufficiently large $m$, the distance of $\partial B_{r_\infty - 2\epsilon}$ from $p_\infty$ with respect to this metric is greater than $(r_\infty - 3\epsilon)$, so that
$$\overline{B}_{r_\infty - 3\epsilon} \subseteq B_{r_\infty - 2\epsilon}.$$

In particular, for all such $m$, this closed ball is compact, and so
$$d(S_m, \hat{e}_m, p_m) \geqslant r_\infty - 3\epsilon.$$

It follows that
$$\mathop{\mathrm{LimInf}}_{m \to \infty} r_m \geqslant r_\infty - 3\epsilon.$$

Since $\epsilon > 0$ is arbitrary, the result follows. $\square$

Define $n : \mathrm{MGL}_{k,c} \to [\sqrt{k}, \infty[$ by
$$n(S, \hat{e}, p) := \|A(p)\|, \tag{0.6}$$

where $A$ here denotes the shape operator of the immersion $(\pi \circ \hat{e})$, and $\|\cdot\|$ denotes its operator norm. We extend this function to one defined over the whole of $\overline{\mathrm{MGL}}_{k,c}$ be setting it equal to infinity over $\mathrm{MT}_c$. This function is trivially continuous and invariant under the action of $\mathrm{Isom}(X_c^3)$, and thus descends to a continuous function over $\mathcal{L}_{k,c}$.

**Lemma 3.3**

For all $B > 0$, there exists $\delta := \delta_1(B, e, k)$ such that, for all $x \in \mathcal{L}_{k,c}$, if $d(x) < \delta$, then $n(x) > B$.

**Proof:** Since $n$ is continuous, $n^{-1}([0, B])$ is closed, and therefore compact. Since $d$ is lower semi-continuous, it attains its minimum $\delta$ over this set. Since this set contains no tubes, $\delta > 0$. It follows that if $d(x) < \delta$, then $x \notin n^{-1}([0, B])$, so that $n(x) > B$, as desired. $\square$

We say that a marked Gauss lift $(S, \hat{e}, p)$ is *umbilic* whenever the immersion $(\pi \circ \hat{e})$ is umbilic at $p$. Let $\mathrm{U}_{k,c}$ denote the set of umbilic points in $\overline{\mathrm{MGL}}_{k,c}$. In particular,

$$\mathrm{U}_{k,c} = n^{-1}(\{\sqrt{k}\}). \tag{0.7}$$

We define two foliations $F_\pm$ of $\overline{\mathrm{MGL}}_{k,c} \setminus \mathrm{U}_{k,c}$ as follows. Let $x := (S, \hat{e}, p) \in \overline{\mathrm{MGL}}_{k,c} \setminus \mathrm{U}_{k,c}$. Let $g$ and $\hat{g}$ denote respectively the metric of $X_c^3$ and the Sasaki metric of $S_1 X_c^3$. Let $M_x$ denote the matrix of $\hat{e}^* \pi^* g$ with respect to $\hat{e}^* \hat{g}$. Since $x$ is non-umbilic, this matrix has distinct eigenvalues. Let $L_{x,-}, L_{x,+} \subseteq T_p S$ denote the respective eigenlines of its lesser and greater eigenvalues. $L_\pm$ define smooth distributions over the set of non-umbilic points in $S$, and the desired foliations are obtained upon integrating these distributions in $S$. When $(S, \hat{e}, p)$ is a marked Gauss lift, $F_-$ and $F_+$ are respectively the lines of greater and lesser curvature of $(\pi \circ \hat{e})$, and when $(S, \hat{e}, p)$ is a tube, $F_-$ and $F_+$ are respectively its longitudinal and transverse curves.

Define $d_+ : \mathrm{MGL}_{k,c} \setminus \mathrm{U}_{k,c} \to ]0, \infty]$ such that, for all $(S, \hat{e}, p) \in \mathrm{MGL}_{k,c} \setminus \mathrm{U}_{k,c}$, $d_+(S, \hat{e}, p)$ is the distance along $F_+$ to $\mathrm{U}_{k,c} \cup \partial_\infty S$. We extend $d_+$ to a function over $\overline{\mathrm{MGL}}_{k,c}$ by setting it equal to zero over $\mathrm{U}_{k,c}$ and infinity over $\mathrm{MT}_c$. This function is trivially invariant under the action of $\mathrm{Isom}(X_c^3)$, and thus descends to a function over $\mathcal{L}_{k,c}$.

**Lemma 3.4**

$d_+$ is lower semi-continuous, that is, for all $x \in \mathcal{L}_{k,c}$

$$\mathop{\mathrm{LimInf}}_{y \to x} d_+(y) \geqslant d_+(x). \tag{0.8}$$

**Proof:** Let $(S_m, \hat{e}_m, p_m)_{m \in \mathbb{N}}$ be a sequence in $\overline{\mathrm{MGL}}_{k,c}$ converging to $(S_\infty, \hat{e}_\infty, p_\infty)$, and let $(\Phi_m)$ be a corresponding sequence of convergence maps. For all $m \in \mathbb{N} \cup \{\infty\}$, denote $d_{+,m} := d_+(S_m, \hat{e}_m, p_m)$ and let $F_{m,+}$ denote the leaf of $F_+$ passing through $p_m$. Choose $\epsilon > 0$, and let $F_{\infty,+}^\epsilon$ denote the closed segment in $F_{\infty,+}$ of length $2(d_{\infty,+} - \epsilon)$ centred on $p_\infty$. Let $\Omega$ be a relatively compact neighbourhood of this segment of $S_\infty$. Let $M$ be such that, for all $m \geqslant M$, the restriction of $\Phi_m$ to $\Omega$ is a diffeomorphism onto its image. Trivially $(\Phi_m^{-1}(F_{m,+}) \cap \Omega)_{m \geqslant M}$ converges in the $C_\mathrm{loc}^\infty$ sense to $F_{\infty,+} \cap \Omega$. In particular, for sufficiently large $m$, $F_{m,+}$ also contains a segment of length $2(d_{\infty,+} - 2\epsilon)$ centred on $p_m$, so that

$$\mathop{\mathrm{LimInf}}_{m \to \infty} d_{m,+} \geqslant d_{\infty,+} - 2\epsilon.$$

Since $\epsilon > 0$ is arbitrary, it follows that

$$\mathop{\mathrm{LimInf}}_{m \to \infty} d_+(S_m, \hat{e}_m, p_m) \geqslant d_+(S_\infty, \hat{e}_\infty, p_\infty),$$

as desired. $\square$

**Lemma 3.5**

For all $B > 0$, there exists $\delta := \delta_2(B, c, k)$ such that, for all $x \in \mathcal{L}_{k,c}$, if $d(x) < \delta$, then $d_+(x) > B$.

**Proof:** Since $d_+$ is lower semi-continuous, $d_+^{-1}([0, B])$ is closed and therefore compact. Since $d$ is lower semi-continuous, it attains a minimum $\delta$ over this set. Since this set contains no tubes, $\delta > 0$. It follows that if $d(x) < \delta$, then $x \notin d_+^{-1}([0, B])$, so that $d_+(x) > B$, as desired. $\square$

Let $x := (S, \hat{e}, p) \in \mathrm{MGL}_{k,c} \setminus \mathrm{U}_{k,c}$. Consider the restriction of $(\pi \circ \hat{e})$ to the leaf of $F_+$ passing through $p$, and define $\kappa_+(S, \hat{e}, p)$ to be the geodesic curvature of this curve at $p$. We extend $\kappa_+$ to a function over the whole of $\overline{\mathrm{MGL}}_{k,c}$ by setting it equal to $-\infty$ over $\mathrm{U}_{k,c}$.

**Lemma 3.6**

$\kappa_+$ *is continuous away from the set of umbilic points.*

**Proof:** Let $(S_m, \hat{e}_m, p_m)$ be a sequence in $\overline{\mathrm{MGL}}_{k,c}$ converging to $(S_\infty, \hat{e}_\infty, p_\infty)$ and let $(\Phi_m)$ be a corresponding sequence of convergence maps. It suffices to consider the case where $(S_\infty, \hat{e}_\infty, p_\infty)$ is a marked tube and $(S_m, \hat{e}_m, p_m)$ is not a marked tube for any finite $m$. Let $\gamma_\infty :\ ]-2\epsilon, 2\epsilon[ \to S_\infty$ be a unit-speed parametrisation of the leaf of $F_+$ passing through $p_\infty$ such that $\gamma_\infty(0) = p_\infty$. By definition $(\pi \circ \hat{e}_\infty \circ \gamma_\infty)$ is a unit-speed parametrised geodesic segment in $X_c^3$. Let $\Omega$ be a relatively compact open subset of $S_\infty$ containing $\gamma_\infty(]-2\epsilon, 2\epsilon[)$. Let $M_1 > 0$ be such that, for all $m \geqslant M_1$, $d_+(S_m, \hat{e}_m, p_m) > \epsilon$. For all such $m$, let $\gamma_m :\ ]-\epsilon, \epsilon[ \to S_m$ be a unit-speed parametrisation of the leaf of $F_+$ passing through $p_m$ such that $\gamma_m(0) = p_m$. Let $M_2 \geqslant M_1$ be such that, for all $m \geqslant M_2$, the restriction of $\Phi_m$ to $\Omega$ is a diffeomorphism onto its image and $\gamma_m(]-\epsilon, \epsilon[) \subseteq \Phi_m(\Omega)$. Trivially, $(\Phi_m^{-1} \circ \gamma_m)_{m \geqslant M_2}$ converges in the $C^\infty_{\mathrm{loc}}$ sense to $\gamma_\infty$. Consequently,

$$\operatorname*{Lim}_{m \to \infty} (\pi \circ \hat{e}_m \circ \gamma_m) = \operatorname*{Lim}_{m \to \infty} (\pi \circ \hat{e}_m \circ \Phi_m) \circ (\Phi_m^{-1} \circ \gamma_m) = \pi \circ \hat{e}_\infty \circ \gamma_\infty,$$

where convergence is here in the $C^\infty_{\mathrm{loc}}$ sense. In particular

$$\operatorname*{Lim}_{m \to \infty} \kappa_+(S_m, \hat{e}_m, p_m) = 0,$$

and the result follows. $\square$

**Lemma 3.7**

*For all $\epsilon > 0$, there exists $\delta := \delta_3(B, c, k)$ such that, for all $x \in \mathcal{L}_{k,c}$, if $d(x) < \delta$, then $\kappa_+(x) < \epsilon$.*

**Proof:** By Lemma 3.3, there exists $\delta_1 > 0$ such that, if $d(x) \leqslant \delta_1$, then $n(x) > \sqrt{k}$, and $x$ is therefore not umbilic. Since $\kappa_+$ is continuous, $d^{-1}([0, \delta_1]) \cap \kappa_+^{-1}([\epsilon, \infty[)$ is closed and thus compact. Since $d$ is continuous, it attains a minimum $\delta_2$ over this set. Since this set contains no tubes, $\delta_2 > 0$. It follows that if $d(x) < \delta_2$, then $x \notin \kappa_+^{-1}([\epsilon, \infty[)$, so that $\kappa_+(x) < \epsilon$, as desired. $\square$

For all $\delta > 0$, we define the open subset $\mathcal{L}_{k,c}^\delta$ by

$$\mathcal{L}_{k,c}^\delta := d^{-1}([0, \delta]). \tag{0.9}$$

By Lemma 3.5, for all $B > 0$, upon choosing $\delta$ sufficiently small, we may suppose that $d_+ \geqslant B$ over $\Omega_\delta$. We then define $d^B(S, \hat{e}, p)$ to be the supremum of $d$ over the closed interval in $F_+$ of length $2B$ centred on $p$. We likewise define $\kappa_+^B(S, \hat{e}, p)$ to be the supremum of $\kappa_+$ over the same closed interval. Proceeding as before, we obtain the following two results.

**Lemma 3.8**

$d^B$ *is lower semi-continuous. In particular, for all $\epsilon > 0$, there exists $\delta > 0$ such that, for all $x \in \mathcal{L}_{k,c}^\delta$, $d^B(x) < \epsilon$.*

**Lemma 3.9**

$\kappa_+^B$ *is continuous. In particular, for all $\epsilon > 0$, there exists $\delta > 0$ such that, for all $x \in \mathcal{L}_{k,c}^\delta$, $\kappa_+^B(x) < \epsilon$.*

**4 - Convexity.** We now construct charts of $S$ in $X_{c+k}^2$. Denote

$$r_0 := \begin{cases} \frac{\pi}{2\sqrt{c+k}} & \text{if } c+k > 0, \text{ and} \\ \infty & \text{otherwise.} \end{cases} \tag{0.10}$$

Note that, when $(c + k) > 0$, $r_0$ is the radius of a hemisphere in $X_{c+k}^2$.

Let $(S, e, p)$ be a marked $k$-surface, and let $g$ denote the metric induced over $S$ by $e$. Recall that $g$ has constant curvature equal to $(c + k)$. Let $\mathrm{Exp}_p : T_p S \to S$ denote the exponential map of $S$ at $p$, and let $\mathrm{Dom}_p$ denote its domain. Recall that $\mathrm{Dom}_p$ is open and star-shaped about $0$. Define $U \subseteq T_p S$ by

$$U := B_{r_0}(0) \cap \mathrm{Dom}_p. \tag{0.11}$$

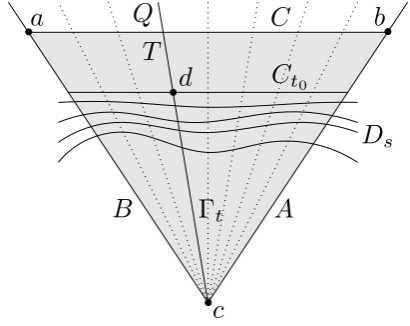

Figure 4.1

Let $\phi : S \to X^2_{c+k}$ be a local isometry. Denote $\hat{\phi} := \phi \circ \mathrm{Exp}_p$. Note that this function restricts to a diffeomorphism from $U$ onto its image $V := \hat{\phi}(U)$. Denote $q := \phi(p)$ and let $\psi : V \to S$ denote the unique left-inverse of $\phi$ such that $\psi(q) = p$.

Recall that a subset of $X^2_{c+k}$ is said to be convex whenever it contains every length-minimising geodesic between any two of its elements.

**Lemma 4.1**

*V is convex.*

**Proof:** Choose $\alpha, \beta \in U$. Denote $a := \hat{\phi}(\alpha)$, $b := \hat{\phi}(\beta)$ and $c := \hat{\phi}(0)$. Consider the closed geodesic triangle $T$ in $X^2_{c+k}$ defined by these three points, as in Figure 4.1. Let $A$, $B$ and $C$ denote respectively the edges lying opposite $a$, $b$ and $c$. We also denote respectively by $A$ and $B$ the extensions of these segments to complete geodesics in $X^2_{c+k}$. Let $Q$ denote the unique quadrant in $M^2_{c+k}$ bounded by $A$ and $B$ and containing $T$. Let $p : [0,1] \to C$ be a parametrisation. For all $t \in [0,1]$, let $\Gamma_t$ denote the unique geodesic in $X^2_{c+k}$ joining $c$ and $p(t)$.

For all $t \in [0,1]$, denote $a_t := \hat{\phi}(t\alpha)$ and $b_t := \hat{\phi}(t\beta)$, and let $C_t$ denote the minimising geodesic segment joining these two points. Denote
$$I := \{t \in [0,1] \mid C_t \not\subseteq V\}.$$

It trivially suffices to show that $I$ is empty. Suppose the contrary. Note that $I$ is closed and therefore contains a least element $t_0 > 0$, say. Let $d$ be a point of $C_{t_0} \setminus V$. Let $\gamma : [0,1] \to T$ denote the unique geodesic joining $c$ and $d$. Note that
$$\lim_{s \to 1}(d \circ \psi \circ \gamma)(s) = 0. \tag{0.12}$$

For all $s$ sufficiently close to 1, let $F_{s,+}$ denote the leaf of $F_+$ passing through $(\psi \circ \gamma)(s)$ and let $D_s$ denote its image under $\phi$. By (0.12) and Lemmas 3.8 and 3.9, $D_s$ converges to a complete geodesic $D_1$, say, in $X^2_{c+k}$.

We first claim that $D_1 \cap T = C_{t_0}$. Indeed, otherwise $D'_1 := D_1 \cap \mathrm{Int}(T_{t_0})$ would be non-trivial. However, by (0.12) and Lemma 3.8, $(d \circ \psi)$ would vanish over this set, which is absurd, and the assertion follows. In particular, there exists $\epsilon > 0$ such that, for all $s \in ]1 - \epsilon, 1]$, and for all $t \in [0,1]$, $D_s$ meets $\Gamma_t$ transversally. We now claim that $(D_s \cap Q)_{s \in ]1-\epsilon,1[}$ foliates an open subset of $T$. Indeed, for all $s \in ]1 - \epsilon, 1[$ and for all $t \in [0,1]$, let $p_{s,t}$ denote the unique point of intersection of $D_s$ and $\Gamma_t$ in $Q$. For all $t \in [0,1]$, $p_{s,t}$ advances strictly monotonically along $\Gamma_t$ as $s$ increases, and it follows that $(D_s \cap Q)_{s \in ]1-\epsilon,1[}$ foliates some open subset $\Omega$, say, of $T$, as asserted. Finally, since $D_1 \cap T = C_{t_0}$, $C_{t_0}$ lies on the boundary of $\Omega$ so that, for all $t \in [0,1]$,
$$\gamma(t) = \lim_{s \to 1} \gamma_{s,t}.$$

In particular, by (0.12) and Lemma 3.8 again,
$$d(\mathrm{Exp}_p(t_0\alpha)) = d(\mathrm{Exp}_p(t_0\beta)) = 0.$$

This is absurd, and the result follows. □

**Lemma 4.2**

$\partial V \cap B_{r_0}(q)$ is a union of geodesics.

**Proof:** Choose $\alpha \in \partial U \cap B_{r_0}(0)$ and define $\gamma : [0,1] \to X^2_{c+k}$ by $\gamma(t) := \hat{\phi}(t\alpha)$. Note that

$$\lim_{t \to 1}(d \circ \psi \circ \gamma)(t) = 0. \tag{0.13}$$

For all $t$ sufficiently close to 1, let $D_t$ denote the image under $\phi$ of the leaf of $F_+$ passing through $(\psi \circ \gamma)(t)$. For all such $t$, since $d$ is non-vanishing over this leaf of $F_+$, $D_t$ can only intersect $\partial V$ along $\partial B_{r_0}(0)$. However, by (0.13) and Lemmas 3.8 and 3.9, $(D_t)_{t \in ]1-\epsilon, 1[}$ converges to a complete geodesic $D_1$, say, in $X^2_{c+k}$. It follows that $\overline{V}$ contains a geodesic passing through $\hat{\phi}(\xi)$. Since $\xi \in \partial U \cap B_{r_0}(0)$ is arbitrary, this completes the proof. $\square$

We now prove our main results.

**Proof of Theorem 1.1:** It trivially suffices to prove the result when $S$ is simply-connected. Note that $\psi$ extends to a homeomorphism from $\overline{V} \cap B_{r_0}(q)$ into an open subset of $\overline{S}$. Since $p \in S$ is arbitrary, this yields an extension of the $X^2_{c+k}$-structure of $S$ to an $X^2_{c+k}$-structure over $\overline{S}$ with geodesic boundary. This extension is trivially unique. Since $e$ is a local isometry with respect to the metric $d$ defined in (0.1), it trivially extends to a continuous function $\overline{e} : \overline{S} \to X^3_c$.

It remains only to show that $\overline{e}$ sends components of $\partial S$ locally isometrically to geodesics in $X^3_c$. Choose $\alpha \in \partial U \cap B_{r_0}(0)$ and define $\gamma : [0,1] \to X^2_{c+k}$ by $\gamma(t) := \hat{\phi}(t\alpha)$. For all $t$ sufficiently close to 1, let $D_t$ denote the image under $\phi$ of the leaf of $F_+$ passing through $(\psi \circ \gamma)(t)$. By Lemmas 3.8 and 3.9, $D_t$ converges to a complete geodesic in $X^2_{c+k}$ which we know lies along $\partial V$. For all $t$ sufficiently close to 1, let $\delta_t : ]-\epsilon,\epsilon[ \to X^2_{c+k}$ be a unit-speed parametrisation of $D_t$ such that $\delta_t(0) = \gamma(t)$. Trivially $(\delta_t)_{t \in ]1-\epsilon,\epsilon[}$ converges to a unit-speed parametrisation of a portion of $\partial V$. On the other hand, since $e$ is a local isometry, by Lemma 3.9, $(\pi \circ \hat{e} \circ \psi \circ \delta_t)_{t \in ]1-\epsilon,1[}$ converges to a unit-speed parametrisation of a geodesic segment in $X^3_c$. It follows that $\overline{e}$ indeed maps $\partial S$ locally isometrically to geodesic segments in $X^3_c$, and this completes the proof. $\square$

**Proof of Corollary 1.2:** Indeed, it suffices to observe that the universal cover of $S$ is isometric to a convex subset of $X^2_{c+k}$. $\square$

**Proof of Theorem 1.3:** We may suppose that $S$ is simply-connected. We first show that $S$ is complete. Indeed, suppose the contrary. We claim that $S$ is isometric to a hemisphere in $X^2_{c+k}$. Indeed, let $\Sigma$ denote the surface obtained by doubling $S$ along its boundary, let $\tilde{\Sigma}$ denote its universal cover, and let $\pi : \tilde{\Sigma} \to \Sigma$ denote the canonical projection. Note that $\pi^{-1}(\partial S)$ is a union of *disjoint* complete geodesics in $\tilde{\Sigma}$. Since $\tilde{\Sigma}$ is complete, it is isometric to $X^2_{c+k}$, from which it follows that $\pi^{-1}(\partial S)$ consists of a single complete geodesic $\Gamma$. Since $S$ identifies with one component of the complement of $\Gamma$, we conclude that $S$ is indeed isometric to a hemisphere in $X^2_{c+k}$, as desired. Finally, let $\overline{e} : \overline{S} \to X^3_c$ denote the continuous extension of $e$. By Theorem 1.1, $\overline{e}$ maps $\partial S$ locally isometrically into a geodesic in $X^3_c$. This is absurd, since $X^3_c$ contains no complete geodesic of length $2\pi/\sqrt{c+k}$, and it follows that $S$ is complete, as asserted.

The result now follows from the classical case where $S$ is complete. We sketch the proof for the reader's convenience. Note first that $S$ is isometric to the sphere $X^2_{c+k}$. Let $I_e$ and $II_e$ denote the first and second fundamental forms of $e$, and note that both define riemannian metrics over $S$. Let $\phi$ denote the Hopf differential - that is, the $(2,0)$-component - of $I_e$ with respect to $II_e$. Since $e$ has constant extrinsic curvature, $\phi$ is holomorphic and therefore vanishes, since $S$ carries no non-trivial holomorphic quadratic differential. $II_e$ is thus a constant multiple of $I_e$, and we conclude by the fundamental theorem of surface theory that $(S,e)$ is a geodesic sphere, as desired. $\square$

**5 - Bibliography.**